\documentclass[leqno]{aomamlt2e}

\numberwithin{equation}{section}
\newenvironment{proof}{\noindent\textit{Proof}: }{\hfill$\square$\vskip 0.5\baselineskip}
\begin{document}

%% title
\title{\flushleft The Second Order Estimate for Fully
Nonlinear Uniformly Elliptic Equations without Concavity
Assumption}
%% author
\author{G.C.Dong, B.J.Bian and Z.C.Guan\thanks{This work was supported by NNSF of
China}}

%\markboth{DRAFT}{}
%% abstract
\begin{abstract}
Investigating for interior regularity of viscosity solutions to
the fully nonlinear elliptic equation
$$F(x,u,\triangledown u,\triangledown ^2 u)=0,$$
we establish the interior $C^{1+1}$ continuity under the
assumptions that $F$ is uniformly elliptic, H$\ddot o$lder
continuous and satisfies the natural structure conditions of
fractional order, but without the concavity assumption of $F$.
These assumptions are weaker and the result is stronger than that
of Caffarelli and Wang[1], Chen[2].
\end{abstract}
%% keywords
\iffalse
\begin{keywords}
 regularity, viscosity solutions, fully nonlinear elliptic
equation, concavity.
\end{keywords}
\fi

%\maketitle

%%content begin

%% Section I
\section{Introduction}

The study of solvability problems for the second order fully
nonlinear uniformly elliptic equation, i.e., the existence,
uniqueness and regularity of solutions with Dirichlet boundary
data can be divided into two classes naturally. The first one is
about classical solvability, i.e., solution $u \in C^2$. The
related results are rich and systematic, See [3], [4], [5] and
therein. All of them is obtained under the assumption of concavity
for the equation with respect to their arguments. The second is
without the concave condition of $F$. In this case, we must seek
for solution in generalized sense. The suitable one is viscosity
solution[6]. But in this case the results are incomplete. Because
the regularity is a key stone for existence and uniqueness, we
prove the interior $C^{1+1}$ continuity without the concavity
assumption of $F$ in this paper and the existence and uniqueness
results will be in next one. Our assumptions are $F(x,z,p,X) \in
C^\beta (\beta
>0)$ of its arguments only, weaker than those in Caffarelli and Wang[1] ($\beta =1$)
and Chen [2]($\beta>1/2$). And the result, $u\in C^{1+1}$ is
stronger than theirs ($u\in C^{1+\alpha},0<\alpha $ and small ).
This paper is organized in five sections: The second is
preliminaries, statements for conditions and the main results. The
third is general comparison principle.
 For solution $u(x)(x\in \Omega, \Omega\subset
\mathbb{R}^n )$, we investigate the general conditions for
$u(x)-u(y)-\Phi(x,y)$ takes maximum in a $\mathbb{R}^{2n}$ domain
$Q\subset \Omega\times \Omega$, where $\Phi(x,y)\in C^2(Q)$. If
these conditions are violate, together with the assumption
$[u(x)-u(y)-\Phi]_{\partial Q}\leq 0$, we have $u(x)-u(y)-\Phi\leq
0$, and the useful estimation $u(x)-u(y)\leq \Phi$ follows. The
fourth is H$\ddot{o}$lder and Lipschitz continuity. The
H$\ddot{o}$lder and Lipschitz estimates are obtained by selecting
suitable $Q$ and $\Phi$. Although these results are not new, we
prove them here for the sake of applications of the results of
section 3.

The last section is $C^{1+1}$ estimate. Having got the interior
Lipschitz estimate for $u(x)$, we can conclude that there exist
Caffarelli points such that in small neighborhood of them $u(x)$
can be separated into linear and second order parts. We estimate
the lipschitz coefficient of $u(x)$ subtracting linear part in
suitable neighborhood of Caffarelli point by selecting suitable
comparison function $\Phi(x,y)$. By this way we get a fundamental
lemma: the lipschitz coefficient is diminishing in certain
constant ratio accompanied with diminishing of radius of spherical
neighborhood in suitably constant ratio.

   The $C^{1+\alpha}$ estimate and $C^{1+1}$ in small part follows
from this fundamental lemma. And the $C^{1+1}$ estimate in the
whole region follows by putting the estimates in small parts
altogether.

%Section II

\section{Preliminaries}
\par Let $\Omega$ be a bounded domain in $\mathbb{R}^n$ with the boundary $ \partial \Omega \in C$.
We consider the regularity problem for solutions of the equation

\begin{equation}\label{eq1}
 F(x,u,\triangledown u,\triangledown ^2 u)=0 \qquad \qquad in \qquad \Omega
\end{equation}
 where $F(x,u,p,X)$ is a function on $\Gamma=\Omega\times
 \mathbb{R}\times \mathbb{R}^n\times \mathbb{S}^n $, $\mathbb{S}^n$ the space of $n\times n$ symmetric
 matrices equipped with usual order. We assume that $F(x,u,p,X)$
 is uniformly elliptic in the following sense
\begin{equation}\label{eq2}
 \lambda Tr(X-Y)\le F(x,z,p,Y)-F(x,z,p,X)\le \Lambda Tr(X-Y)\end{equation}
 for all $(x,z,p,X)\in \Gamma$, $Y\in \mathbb{S}^n$ and $X\ge Y$, where
 $\lambda$ and $\Lambda$ are positive constants with $\lambda \le
 \Lambda$. Assume that $F$ is monotone with respect to z

 \begin{equation}\label{eq3}
 F(x,z,p,X)-F(x,w,p,X)\le 0
\end{equation}
$\forall z\le w$. Furthermore, we suppose that there exist
positive constants $\mu$ and $\beta$ with $\beta \le 1$ such that
\begin{equation}\label{eq4}
\begin{aligned}
  |F(x,z,p,X)-F(y,w,q,X)|\qquad\qquad\\
  \le \mu [|x-y|(1+|p|+|q|)+\frac {|p-q|}{1+|p|+|q|}]^\beta\\
       \cdot (1+|p|^2+|q|^2+||X||)\qquad \qquad\qquad
 \end{aligned}
\end{equation}
$\forall$ $(x,z,p,X)\in \Gamma$ and $(y,w,q,X)\in \Gamma $.
(\ref{eq4}) is called natural structural condition for $F$ of
fractional order (order $\beta$, where $0<\beta\leq 1$).

Now we give the definition of viscosity solutions and main
theorem.

\par {\bf Definition 1.}
  Let $u$ be an upper semi-continuous (resp. lower
semi-continuous) function in $ \Omega $. $u$ is said to be a
viscosity subsolution (resp. supersolution) of ($\ref{eq1}$) if
for all $\phi \in C^2(\Omega)$ the following inequality
$$F(x_0,u(x_0),\triangledown \phi (x_0),\triangledown ^2 \phi (x_0))\le 0$$
$$(resp.\qquad F(x,u(x_0),\triangledown \phi(x_0),\triangledown ^2 \phi(x_0))\ge 0)$$
holds at each local maximum (resp. minimum) point $x_0\in \Omega$
of $u-\phi$.

\par {\bf Definition 2.}  $u\in C(\Omega)$ is said to be a viscosity solution
of ($\ref{eq1}$) if $u$ is both a viscosity subsolution and a
supersolution.

\par {\bf Theorem 2.1.} Assume $F(x,z,p,X)$ satisfies the conditions $($\ref{eq2}$)-($\ref{eq4}$)$ and $u$ is
a viscosity solution of ($\ref{eq1}$), then $ u \in
C^{1+1}(\Omega)$.

In this paper, we don't study any boundary value problem for
(\ref{eq1}). About the existence and uniqueness for solution of
boundary value problems, we shall discuss them in a next paper.

\section{General Comparison Principle}
\par For suitably selected regular non-negative function $\Phi (x,y)$ and
a bounded domain $Q \subset \mathbb{R}^{2n}$, we suppose that
$u(x)-u(y)-\Phi (x,y)$takes a positive maximum value at an
interior point $(\bar x,\bar y)\in Q$, then
$$  u(x)-u(y)-\Phi (x,y) \le u(\bar x)-u(\bar y)-\Phi (\bar x,\bar
y)> 0$$

For simplicity, we omit the upper bar on $x,y$ in the following,
i.e., write $(x,y)=(\bar x,\bar y)$, thus
\begin{equation}\label{eq5}
 u(x)-u(y)> \Phi (x,y)\ge 0.
\end{equation}
In particular we have, $x\neq y$. From [6], there exist $X,Y \in
\mathbb{S}^n$ such that
\begin{equation}\label{eq6}
 F(x,u(x),\Phi_x(x,y),X)-F(y,u(y),- \Phi_y(x,y),-Y)\le0
\end{equation}
and
\begin{equation}\label{eq7}
\left(
\begin{array}{cc}
 X&0\\0&Y
 \end{array}
 \right)
 \le
\left(
 \begin{array}{cc}
 {\Phi _{xx}}&{\Phi _{xy}}\\ {\Phi _{yx}}&{\Phi _{yy}}
\end{array}
\right)
\end{equation}
From above inequality, we have
\begin{equation}\label{eq8}
\left(
\begin{array}{cc}
 X+Y&X-Y\\X-Y&X+Y
 \end{array}
 \right)
 \le
\left(
 \begin{array}{cc}
 Z_1&Z\\Z&Z_2
\end{array}
\right),
\end{equation}
where
\begin{equation}\label{eq9}
\begin{aligned}
Z_1=(\frac{\partial}{\partial x}+\frac{\partial}{\partial y})^2\Phi ,\\
Z_2=(\frac{\partial}{\partial x}-\frac{\partial}{\partial y})^2\Phi ,\\
Z=(\frac{\partial}{\partial x}-\frac{\partial}{\partial y})(\frac
{\partial}{\partial x}+\frac{\partial}{\partial y})\Phi .
\end{aligned}
\end{equation} A consequence of (\ref{eq8}) are
\begin{equation}\label{eq10}
X+Y-Z_1\le 0,\qquad X+Y-Z_2\le 0.
\end{equation}
$\forall$ $\sigma\in \mathbb{R}$ and $\xi \in \mathbb{R}^n $,
multiplying (\ref{eq8}) from left and right by $(\sigma \xi ,\xi
)$ and $(\sigma \xi ,\xi )^T $ respectively, we get
%\label{eq11}
$$
\begin{aligned}
\sigma ^2<(X+Y-Z_1)\xi ,\xi >+2\sigma <(X-Y-Z)\xi
,\xi> \\
 +<(X+Y-Z_2)\xi ,\xi>\le 0.\qquad \qquad \end{aligned}$$
 Hence, for $\xi \in \mathbb{R}^n $ with $ |\xi |=1 $,
%\label{eq12}
$$
\begin{aligned}
<(X-Y-Z)\xi,\xi>^2\le \qquad \qquad  \\
 <(X+Y-Z_1)\xi ,\xi ><(X+Y-Z_2)\xi ,\xi>.
 \end{aligned}$$
 i.e.,
$$
%\label{eq13}
\begin{aligned}
||X-Y-Z||^2\le \qquad \qquad  \\
 ||X+Y-Z_1|| ||X+Y-Z_2||\le \\
 ||X+Y-Z_1||^2+||Z_2-Z_1||||X+Y-Z_1||.\end{aligned}$$
This inequality implies that
$$||X-Y-Z||$$
$$\le ||X+Y-Z_1||+||Z_2-Z_1||^{1/2}||X+Y-Z_1||^{1/2}$$
%\begin{equation}\label{eq14}
$$
 \le C_1[|Tr(X+Y-Z_1)|+||Z_2-Z_1||^{1/2}|Tr(X+Y-Z_1)|^{1/2}].$$
 Thus,
$$||X||+||Y|| $$
$$\le ||X+Y||+||X-Y||$$
$$\le ||Z_1||+||Z||+||X+Y-Z_1||+||X-Y-Z||$$
$$\le C_1[||Z_1||+||Z||+2|Tr(X+Y-Z_1)| $$
%\begin{equation}\label{eq15}
$$
 +||Z_2-Z_1||^{1/2}|Tr(X+Y-Z_1)|^{1/2}].
$$
 Taking a positive parameter $\omega $ to be determined later. From the above inequality, we have
\begin{equation}\label{eq16}\begin{aligned}
\omega (||X||+||Y||) \qquad \qquad  \qquad \qquad \\
  \le C_1\omega[||Z_1||+||Z||+2|Tr(X+Y-Z_1)|] \\
 +\frac {\lambda}{2}|Tr(X+Y-Z_1)|+\frac {2C_1^2\omega
 ^2}{\lambda}||Z_2-Z_1||.
 \end{aligned}\end{equation}
 On the other hand, by ($\ref{eq2}$) and ($\ref{eq10}$), we have
$$\lambda |Tr(X+Y-Z_1)|=\lambda Tr(Z_1-X-Y)$$
 $$\le F(x,u(x),\Phi _x,X)-F(x,u(x),\Phi _x,Z_1-Y).$$

Applying ($\ref{eq6}$) and ($\ref{eq3}$), we have
$$ F(x,u(x),\Phi _x,X)\le F(y,u(y),-\Phi _y,-Y)$$
$$\le  F(y,u(x),-\Phi _y,-Y).$$
Separating $Z_1$ into $Z_1^++Z_1^-$ where $Z_1^+\ge 0$, $Z_1^-\le
0$ and by (\ref{eq4}), we obtain
$$ -F(x,u(x),\Phi _x,-Y+Z_1)$$
$$\le -F(x,u(x),\Phi _x,-Y) +\Lambda TrZ_1^+-\lambda TrZ_1^-$$
Combining the above inequalities and applying (\ref{eq4}), we have
$$\lambda |Tr(X+Y-Z_1)| $$
$$\le C_2\{||Z_1||+ [(1+A)|x-y|+\frac {|Z_0|}{1+\Lambda}]^\beta $$

\begin{equation}\label{eq17}
 \cdot(1+A^2+||X||+||Y||)\},
 \end{equation}
 where
\begin{equation}\label{eq18}\begin{aligned}
Z_0=\Phi_x+\Phi_y=(\frac {\partial}{\partial x}+\frac {\partial}{\partial y})\Phi ,\\
 A=|Z_0|+|(\frac {\partial}{\partial x}-\frac {\partial}{\partial
 y})\Phi|.
 \end{aligned}\end{equation}

Taking the parameter $\omega$ to be
\begin{equation}\label{eq19}
\omega=C_2[(1+A)|x-y|+\frac {|Z_0|}{1+\Lambda}]^\beta.
\end{equation}
Combining (\ref{eq16}), (\ref{eq17}) and (\ref{eq19}), we have
\begin{equation}\label{eq20}\begin{aligned}
(\frac \lambda 2 -2C_1\omega)|Tr(X+Y-Z_1)|\qquad \qquad \\
\le C[||Z_1||+(1+A^2+||Z||)\omega+||Z_2-Z_1||\omega ^2].
\end{aligned}\end{equation}
As $\omega$ is small, $\omega \ll 1$, (\ref{eq20}) represent an
upper estimate for $|Tr(X+Y-Z_1)|$. We still need a lower estimate
for $|Tr(X+Y-Z_1)|$.

Let $P$ be a $n\times n$ diagonal matrix with $0<P\le I$, where
$I$ is unit matrix. Since $Tr$ is invariant under coordinate
rotation, denote $S$ be coordinate rotation matrix, which is
symmetric and satisfies $S^2=I$, then
$$-Tr[SPS(Z_2-Z_1)]\qquad \qquad$$
$$=-Tr[S^2PS(Z_2-Z_1)S]=-Tr[PS(Z_2-Z_1)S]$$
$$=Tr[PS(Z_1-X-Y)S]-Tr[PS(Z_2-X-Y)S]$$
Applying (\ref{eq10}), we get $$S(Z_1-X-Y)S\ge 0,\ S(Z_2-X-Y)S\ge
0$$

Since $P$ is diagonal and $0<P\le I$, we have $$PS(Z_1-X-Y)S \le
S(Z_1-X-Y)S,\ PS(Z_2-X-Y)S\ge 0.$$ Hence
$$-Tr[SPS(Z_2-Z_1)]\le Tr[S(Z_1-X-Y)S]$$
$$=Tr(Z_1-X-Y)=|Tr(X+Y-Z_1)|$$
Setting
$$ \Upsilon=\frac {(x-y)\otimes (x-y)}{|x-y|^2}$$ it satisfies $0\le \Upsilon\le I$,
 and selecting $S$ such that $SPS=\frac{1}{1+\varepsilon}(\Upsilon+\varepsilon I)$and let $\varepsilon\to 0$,we have
\begin{equation}\label{eq21}
-Tr[\Upsilon(Z_2-Z_1)]\le |Tr(X+Y-Z_1)|.
\end{equation}
(\ref{eq21}) is a lower estimate for $Tr(X+Y-Z_1)$ which are
needed. Combining (\ref{eq20}) and (\ref{eq21}), we have
\begin{equation}\label{eq22}\begin{aligned}
-(\frac \lambda 2-2C_1\omega)Tr(\Upsilon Z_2) \qquad\qquad \\
\le C[||Z_1||+(1+A^2+||Z||)\omega +||Z_2||\omega^2].
\end{aligned}\end{equation}
which is a necessary condition for $u(x)-u(y)-\Phi (x,y)$ takes
positive maximum value in $Q$. If (\ref{eq22}) is violate, i.e.
\begin{equation}\label{eq23}\begin{aligned}
-(\frac \lambda 2-2C_1\omega)Tr(\Upsilon Z_2) \qquad\qquad \\
> C[||Z_1||+(1+A^2+||Z||)\omega +||Z_2||\omega^2],
\end{aligned}\end{equation}
then  $u(x)-u(y)-\Phi (x,y)$ cannot take positive maximum value in
$Q$, hence $\forall \ (x,y)\in Q$, we have
\begin{equation}\label{eq24}
u(x)-u(y)-\Phi (x,y)\le [u(x)-u(y)-\Phi (x,y)]|_{\partial Q}
\end{equation}

The above all are general discuss.

\section{H$\ddot o$lder and Lipschitz continuity}
\label{sec:Lipschitzcontinuity}

{\bf Theorem 4.1} Let $u$ be a viscosity solution to (\ref{eq1})
in $\Omega$ and suppose that (\ref{eq2})-(\ref{eq4}) hold, then
there exists a constant $\alpha \in (0,1)$, such that $u$ is
H$\ddot o$lder continuous with exponent $\alpha$ in $\Omega $ and
the following estimate holds $\forall$ $x,y\in \Omega$,
\begin{equation}\label{eq25}
|u(x)-u(y)|\le \frac C{d^\alpha}|x-y|^\alpha,
\end{equation}
where $d= min[d(x,\partial \Omega),d(y,\partial \Omega)]$ and $C$
is a constant depending only on $n,\lambda,\Lambda,\mu,\beta$ and
$sup_{x\in \Omega}|u(x)|$.

%{\bf Proof.}
\begin{proof} Let $x_0\in \Omega$ and $R>0$ such
that $B(x_0,2R)\subset \Omega$. Without loss of generality,
suppose $x_0$ is the origin and $sup_{B_{2R}}|u(x)|=1$. Let
$\alpha<1$ and $K>1$ be two constants to be chosen and taking the
domain
$$ Q=\{(x,y)\in \Omega ||x|^2+|y|^2<R^2\}$$
and the function \begin{equation}\label{eq226}\Phi (x,y)=\frac
2{R^2}(|x+y|^2+|x-y|^2)+K\frac
{|x-y|^\alpha}{R^\alpha}\end{equation} and setting
$$ w=u(x)-u(y)-\Phi (x,y),$$ then we have $w\le0$ on $\partial Q$.
We claim that $w\le 0$ in $Q$. If it is not true, there exists a
positive maximum value of $w$ in $Q$  at $(x,y)$ and by
(\ref{eq5}) and (\ref{eq6}), we have
%\begin{equation}\label{eq26}
$$
|x-y|\le (\frac {|u(x)-u(y)|}{K})^{1/\alpha} R\le (\frac 2
K)^{1/\alpha}R,$$ or
\begin{equation}\label{eq26}
K\left( \frac{|x-y|}{R}\right)^\alpha\le 2.
\end{equation}

It is easy to see
%\begin{equation}\label{eq27}
$$|Z_0|\le CR^{-1},$$
$$
\begin{aligned}
\frac {\alpha K}{R^\alpha}|x-y|^{\alpha-1}\le A \le
C(R^{-1}+\alpha K \frac {|x-y|^{\alpha -1}}{R^\alpha })
\end{aligned}$$
where $Z_0$ and $A$ are defined by (\ref{eq18}). From (\ref{eq19})
and (\ref{eq26}), we have
%\begin{equation}\label{eq28}
$$
\begin{aligned}
\omega\le C[\frac {|x-y|}{R} +\alpha K (\frac {|x-y|} R)^\alpha
+(\alpha K)^{-1}(\frac {|x-y|} R)^{1-\alpha}]^\beta\\
\end{aligned}$$
$$
\le C(\alpha +\alpha^{-1}K^{-1/\alpha})^\beta=o(1),$$ as taking
$\alpha$ small first and then taking $K$ large. By the definition
of (\ref{eq9}), we get
%\begin{equation}\label{eq29}
$$\begin{aligned}
||Z_1||\le CR^{-2}, Z=0,\\
||Z_2||\le C R^{-2}+C\alpha K R^{-\alpha}|x-y|^{\alpha-2}.
\end{aligned}$$
Hence, for our $\Phi$, by applying (\ref{eq26}), (\ref{eq22})
becomes
\begin{equation}\label{eq30}\begin{aligned}
\frac {4\alpha(1-\alpha)K}{R^\alpha}|x-y|^{\alpha-2}+O(R^{-2})\\
\le C[R^{-2}+(\alpha K R^{-\alpha}|x-y|^{\alpha-1})^2\omega
+\alpha KR^{-\alpha}|x-y|^{\alpha-2}\omega^2]\\
\le C[R^{-2}+\alpha^2 K R^{-\alpha}|x-y|^{\alpha-2}\omega +\alpha
KR^{-\alpha}|x-y|^{\alpha-2}\omega^2].
\end{aligned}\end{equation}
(\ref{eq30}) can not hold when we take $\alpha$ small first and
then take $K$ large. Thus $u(x)-u(y)-\Phi$ cannot takes positive
maximum value in $Q$. Combining estimate for
$[u(y)-u(x)-\Phi]_{\partial Q}$, we have
$$|u(x)-u(y)|\le \Phi=\frac 2{R^2}(|x+y|^2+|x-y|^2)+K\frac
{|x-y|^\alpha}{R^\alpha}.$$
 Especially,
$$|u(x)-u(0)|\le \frac 4{R^2}|x|^2+K\frac
{|x|^\alpha}{R^\alpha}\le (4+K)\frac {|x|^\alpha}{R^\alpha},$$ By
coordinate translation, $\forall$ $y$ such that $|y|<R$, we have
$$|u(x+y)-u(y)|\le (4+K)\frac {|x|^\alpha}{R^\alpha}.$$
The theorem follows by substituting $x+y$ for $x$.
\end{proof}

{\bf Theorem 4.2} Let $u$ be a viscosity solution of (\ref{eq1})
in $\Omega$ and suppose that $(\ref{eq2})-(\ref{eq4})$ hold, then
u is locally Lipschitz continuous and satisfies the following
estimate
\begin{equation}\label{eq31}
|u(x)-u(y)|\le \frac C d |x-y|.\end{equation}

%{\bf Proof.}
\begin{proof}
 Let $x_0\in \Omega$ and $R>0$ such that
$B(x_0,2R)\subset \Omega$. Without loss of generality, suppose
$x_0$ is the origin and $sup_{B_{2R}}|u(x)|=1$. Consider the
function
\begin{equation}\label{eq32}\begin{aligned}
\Phi(x,y)=\frac 2{R^2}(|x+y|^2+|x-y|^2)+4[1-\frac 1
4(\frac{K|x-y|}{R})^\gamma]\frac{K|x-y|}R\end{aligned}\end{equation}
in the domain
$$Q=\{(x,y)\in \Omega||x|^2+|y|^2<R^2, K|x-y|<R\}.$$
where $\gamma(\le \alpha\beta)$ and $K(\ge 1)$ are positive
constants to be chosen. Obviously, we have
$w=u(x)-u(y)-\Phi(x,y)\le 0$ on $\partial Q$. If $w$ takes
positive maximum in $Q$ at (x,y), by $(\ref{eq5}),(\ref{eq25})$
and $(\ref{eq32})$, we have
$$\frac {2|x+y|^2}{R^2}\le \Phi\le u(x)-u(y)\le \frac
{C|x-y|^\alpha}{R^\alpha},$$ hence
$$|x+y|\le CR(\frac{|x-y|}R)^{\alpha/2}$$
and
$$K\frac {|x-y|}{R}\le \Phi\le u(x)-u(y)\le \frac
{C|x-y|^\alpha}{R^\alpha},$$
$$K(\frac {|x-y|}{R})^{1-\alpha}\le
C,\ K^{\frac 1{1-\alpha}}\frac {|x-y|}{R}\le C,\ K\frac
{|x-y|}R\le CK^\frac {-\alpha}{1-\alpha}.$$ Thus
$$K(\frac {|x-y|}{R})\le C K^{-\frac \alpha{2(1-\alpha)}}(\frac
{|x-y|}R)^{\frac \alpha 2},$$ and the corresponding  $Z_0$, $A$
and $\omega $  defined by (\ref{eq18}) and (\ref{eq19}) have the
following estimates.
$$|Z_0|=O(R^{-1}(\frac {|x-y|}R)^\frac \alpha 2),\ \frac K R\le
A=O(\frac K R),$$
$$\omega=O([(1+A)|x-y|+\frac {|Z_0|}{1+A}]^\beta)$$
$$=O(K^{-\frac {\alpha\beta}{2(1-\alpha)}}(\frac {|x-y|}R)^\frac
{\alpha\beta}2)=o(1)$$ as $K$ is large. And
$$||Z_1||=O(R^{-2}),\ Z=0,\ ||Z_2||=O(\frac K{R|x-y|}),$$
$$-Tr(\Upsilon Z_2)=O(R^{-2})+(1+\gamma)\gamma\frac {K^2}{R^2}(\frac
{K|x-y|}R )^{\gamma-1},$$ where $Z,Z_1,Z_2$ are defined by
(\ref{eq9}).

Substituting the above estimates into (\ref{eq22}), we obtain
\begin{equation}\label{eq33}\begin{aligned}
 [\frac \lambda 2+o(1)][O(R^{-2})+(1+\gamma)\gamma\frac
{K^2}{R^2}(\frac {K|x-y|}R )^{\gamma-1}]\\
\le C[O(R^{-2})+K^{-\frac {\alpha\beta}{2(1-\alpha)}}(\frac
{|x-y|}R)^{\frac {\alpha\beta} 2} \frac {K^2}{R^2}\\
+K^{-\frac {\alpha\beta}{1-\alpha}}(\frac
{|x-y|}R)^{\alpha\beta}\frac K{R|x-y|}.
\end{aligned}
\end{equation}
Since $\frac {K|x-y|}R\le C$, then
$(\frac {K|x-y|}R)^{\gamma -1}\ge C$. If we take
$\gamma=\alpha\beta$ and $K$ large enough, (\ref{eq33}) cannot
hold. This means
$$|u(x)-u(y)|\le \Phi (x,y).$$
In particular, we have
$$|u(x)-u(0)|\le \Phi (x,0)$$
$$\le \frac {4|x|^2}{R^2}+\frac {4K|x|}R\le 4(1+K)\frac{|x|}R.$$
Applying coordinate translation, $\forall$ $y$ such that $|y|<K$
we have
$$|u(x+y)-u(y)|\le 4(1+K)\frac{|x|}R.$$
 The theorem is proved by
substituting $x+y$ for $x$.
\end{proof}

\begin{center}
\section{$C^{1+1}$ interior estimate}
\end{center}
\par {\bf Definition 3.} $u$ is said in $C^{1+1}(\Omega)$,
if $Du(x)$ exists $\forall$ $x\in \Omega$ and moreover, $Du(x)$
satisfies Lipschitz condition for all closed subset $\tilde \Omega
\subset \subset \Omega$, and, for all line segment
$\overline{xy}\in\tilde\Omega$,
\begin{equation}\label{eq34}
 |Du(x)-Du(y)|
\le C|x-y|,
\end{equation}
where $C$ depends only on $n,\lambda,\Lambda,\mu,\beta$ and
$dist(\tilde \Omega,\partial\Omega)$.

If (\ref{eq34}) is replaced by H$\ddot o$lder condition and
\begin{equation}\label{eq35}
 |Du(x)-Du(y)|
\le C|x-y|^\alpha (0<\alpha<1),
\end{equation}
we call $u\in C^{1+\alpha}(\Omega).$

The following is a covering theorem.
\par {\bf Theorem 5.1} For any given sufficiently small positive
constant $\theta$ and $\forall\ B(x_0,R_0)\subset \Omega$ with
$R_0\le\theta$, $u\in C^{1+1}(B(x_0,\xi R_0))$ where $\xi$ is a
constant $\xi\in(0,1)$, then $u\in C^{1+1}(\Omega).$

\begin{proof}
%\par {\bf Proof}
For any $\tilde \Omega \subset \subset \Omega$ and all $x,y\in
\tilde\Omega$ with $\bar{xy}\in \tilde \Omega $, we cover
$\bar{xy}$ by a finite set of spheres $\{B(z_i,\xi R_i)\}_{0\le
i\le m-1}\subset \tilde \Omega$, with
$B(z_i,R_i)\subset\tilde\Omega$, where $R_i\le \theta$ with
$x_0=x,x_m=y$ and $\bar{x_ix_{i+1}}\in B(z_i,\xi R_i))(0\le i\le
m-1)$. Hence we have
$$|Du(x)-Du(y)|=|\sum_{i=0}^{m-1}[Du(x_i)-Du(x_{i+1})]|$$
$$\le \sum _{i=0}^{m-1}|Du(x_i)-Du(x_{i+1})|\le C\sum_{i=0}^{m-1}|x_i-x_{i+1}|= C|x-y|.$$
The theorem is proved.
\end{proof}

Fix a $B(x_0,R_0)$ such that $\overline{B(x_0,R_0)}\subset\Omega$
and by section \ref{sec:Lipschitzcontinuity}, we see that there
exists a constant $M_1$ such that $\forall$ $x,y\in B(x_0,R_0)$
\begin{equation}\label{eq36}
 |u(x)-u(y)|
\le M_1|x-y|.\end{equation} Without loss of generality we assume
$\theta \le M_1$,otherwise substituting $M_1$ by
max$(\theta,M_1)$. We consider a function $v(y)$ in $B(0,1)$ as
follows
\begin{equation}\begin{aligned}
 v(y)=\frac {u(x_0+R_0y)-u(x_0)}{\widetilde{M}},\qquad\qquad\\
\widetilde{M}=osc_{x\in
 B(x_0,R_0)}u(x)\qquad\qquad\\
 +R^2_0[1+sup_{|p|\le M_1}sup_{x\in
 B(x_0,R_0)}|F(x,u(x),p,0)|]
.\end{aligned}\end{equation} By the definition of viscosity
solution it is easy to see in $B(0,1)$,

\begin{equation}\label{eq37} |v|\le 1,\ v\in S(f),\ |f|\le 1,\end{equation}
where $S(f)$ denotes the class of viscosity solutions to elliptic
equation related to Pucci's extremal operator(see [7]).

It is well known that (see[8]) the set $Y_E$ of points$\{y_1\}$ in
$B(0,1)$ satisfying the following inequality
\begin{equation}\label{eq38}
 |v(y)-v(y_1)-<a,y-y_1>|
\le E|y-y_1|^2.\end{equation} $\forall$ $y\in B(0,1)$ has density
$\frac{|Y_E|}{|B(0,1)|}$greater than $1-E^{-\Gamma}$, where the
constant $E$ is suitably large, i.e., $E$ is bounded below by a
constant $E_0$ depending only on $n,\lambda,\Lambda,\mu$ and
unbounded above. And  $\Gamma$ satisfies  $0<\Gamma<1$ and depends
only on $n,\lambda,\Lambda,\mu$ as well. All $y_1$ $\in$ $Y_E$ are
called Caffarelli points in $B(0,1)$.

In the following we always fix a large $E$ for our study.
 Taking any point $y_1\in Y_E$ and scaling back as $x_1$, we have $x_1\in
B(x_0,R_0)$ such that $\forall$ $x\in B(x_0,R_0)$
\begin{equation}\label{eq39}\begin{aligned}
 |u(x)-u(x_1)-<a,x-x_1>|\qquad\qquad\\
\le \frac E{R_0^2} M|x-x_1|^2\le \frac {C_1EM_1}{R_0}|x-x_1|^2
 ,\end{aligned}\end{equation}
where
\begin{equation}\label{eq40}\begin{aligned}
 C_1=3+sup_{|p|\le M_1}sup_{x\in
 B(x_0,R_0)}|F(x,u(x),p,0)|.\end{aligned}\end{equation}
The inequality (\ref{eq39}) is called Caffarelli expansion of
$u(x)$ in $B(x_0,R_0)$ and $x_1$ is called a Caffarelli point in
$B(x_0,R_0)$, and (\ref{eq39}) is valid under the restriction
$\theta\leq M_1$, hence $R_0\le \theta\le M_1$, and the constant
vector $a$ can be determined by dividing (\ref{eq39}) by
$|x-x_1|$, and then let $x\rightarrow x_1$ in any fixed direction
$l=\frac {x-x_1}{|x-x_1|}.$ We have
$$<Du(x_1),l>-<a,l>=0,$$
or $a=Du(x_1)$. This means that at any Caffarelli point $x_1$,
$Du(x_1)$ exists. Moreover,
\begin{equation}\label{eq41}
|Du(x_1)|\le \lim_{y\rightarrow x_1}\frac
{|u(x_1)-u(y)|}{|x_1-y|}\le M_1.\end{equation}

Since Caffarelli point of $B(x_0,R_0)$ is always a Caffarelli
point of $B(y_0,S_0)$, where $B(y_0,S_0)\subset B(x_0,R_0)$. And
the inverse is not true. So that talk about Caffarelli point, we
must point out its related sphere simultaneously.
   Let $x_1$ be a Caffarelli point of $B(x_0,R_0)$ such that $B(x_1,\sqrt{2}R_1)\subset B(x_0,R_0)$ and
assume that
\begin{equation}\label{eq42}
|u(x)-u(y)-<Du(x_1),x-y>|\le K|x-y|.\end{equation} $\forall x,y
\in B(x_1,\sqrt{2}R_1).$  Applying (\ref{eq36}) and (\ref{eq41})
we have the Lipschitz coefficient $K\le 2M_1$.

We estimate the Lipschitz coefficient in the following region of
$\mathbf{R}^{2n}$.
\begin{equation}\label{eq43}
Q:\{x,y|J\equiv
{(\frac{|x+y-2x_1|}{2R_1})}^2+{(\frac{|x-y|}{2\epsilon
R_1})}^{\frac{2\sigma}{\ln E}}+e^{-2m(1+2\delta)\sigma} < 1\}
\end{equation}
where $\epsilon, \sigma$ are small positive constants and only
$\epsilon$ is depending on $E$, $\delta$ is a constant, $\delta\in
[\frac 1 2, \frac 3 2], m(>1)$ is a constant independent of $E$.
It is easy to show that
$$Q\subset\{x,y|{(\frac{|x+y-2x_1|}{2R_1})}^2+{(\frac{|x-y|}{2R_1})}^2
< 1\}$$ $$\subset \{x||x-x_1|<\sqrt{2}R_1\}\times
\{y||y-x_1|<\sqrt{2}R_1\}$$ $$\subset \{x||x-x_1|<R_0\}\times
\{y||y-x_1|<R_0\}$$ Take the comparison function $\Phi(x,y)$ to be
\begin{equation}\label{eq44}
\Phi(x,y)=KJ^{\delta/2}[(1+\varphi)|x-y|-\Psi(|x-y|)],
\end{equation} where $\Psi$ is a $C^2$ function to be determined
which satisfies
\begin{equation}\label{eq45}
|\Psi'(|x-y|)| \le \varphi, \qquad \Psi''(|x-y|)=O(\frac1{|x-y|}).
\end{equation}
We investigate the conditions for
\begin{equation}\label{eq46}
|u(x)-u(y)-<Du(x_1),x-y>| \le \Phi(x,y),
\end{equation}
$\forall x,y \in Q$

On $\partial Q$, applying (\ref{eq42}),(\ref{eq44}) and
(\ref{eq45}), we see that (\ref{eq46}) is satisfied.

In $Q$ we have $$\begin{aligned} |Z_0|=|(\frac{\partial}{\partial
x}+\frac{\partial}{\partial y})[\mp(Du(x_1))+\Phi]|\qquad
\qquad\\=4K\delta
J^{\frac{\delta}{2}-1}\frac{|x+y-2x_1|}{R_1^2}[(1+\varphi)|x-y|-\Psi(|x-y|)]\\
\le 4(1+\varphi)K\delta J^{\frac{\delta-1}2}\frac{|x-y|}{R_1}\qquad \qquad \\
=O(K\epsilon{(\frac{|x-y|}{2\epsilon R_1})}^{1-\frac {1}{2}
\frac{\sigma}{\ln E}}).\qquad \qquad \\
A=|Z_0|+|(\frac{\partial}{\partial x}-\frac{\partial}{\partial
x})[\mp<Du(x_1),x-y>+\Phi]|=O(1).\\
\omega=C_2{[(1+A)|x-y|+\frac{|Z_0|}{1+A}]}^\beta\qquad\qquad\\
=O({[(\epsilon R_1+K\epsilon){(\frac{|x-y|}{2\epsilon
R_1})}^{1-\frac {1}{2} \frac{\sigma}{\ln E}}]}^\beta)\qquad\qquad\\
=O({(K\epsilon)}^\beta {(\frac{|x-y|}{2\epsilon
R_1})}^{\beta[1-\frac {1}{2} \frac{\sigma}{\ln E}]})=o(1).\qquad
\end{aligned}$$
When $\epsilon$ is small and we restrict $R_1$ by
\begin{equation}\label{eq47}
R_1 \le K.
\end{equation}
It is easy to calculate
$$\|Z_1\|=O(\frac{K}{R_1^2}J^{\delta/2-1}|x-y|),$$
$$\|Z\|=O(\frac{K}{R_1}J^{\frac{\delta-1}{2}}),$$
$$\|Z_2\|=O(KJ^{\delta/2}\frac 1 {|x-y|}),$$
by using (\ref{eq45}).
\begin{equation}\label{eq48}
-Tr(\Upsilon Z_2)\geq 4KJ^{\delta/2}\Psi''(|x-y|)+
O(KJ^{\delta/2-1}\frac{\sigma}{\ln E}\frac 1 {2\epsilon R_1}
{(\frac{|x-y|}{2\epsilon R_1})}^{\frac{2\sigma}{\ln E}-1}).
\end{equation}

The validity of (\ref{eq48}) is due to
$$ \frac\delta 2(1-\frac\delta 2)J^{\frac\delta 2-2}Tr\{\Upsilon[(\frac{\partial}{\partial
x}-\frac{\partial}{\partial y}){(\frac{|x-y|}{2\epsilon
R_1})}^{\frac{2\sigma}{\ln E}}]$$
$$\otimes[(\frac{\partial}{\partial x}-\frac{\partial}{\partial
y}){(\frac{|x-y|}{2\epsilon R_1})}^{\frac{2\sigma}{\ln
E}}]\}[(1+\varphi)|x-y|-\Psi(|x-y|)]\geq 0, $$ and
$$(\frac{\partial}{\partial
x}-\frac{\partial}{\partial y})\Psi=\Psi'\frac{2(x-y)}{|x-y|},$$
$${(\frac{\partial}{\partial
x}-\frac{\partial}{\partial
y})}^2\Psi=2\Psi'(\frac{\partial}{\partial
x}-\frac{\partial}{\partial y})\frac{x-y}{|x-y|}+
4\Psi''\frac{(x-y)\otimes(x-y)}{{|x-y|}^2}$$
$$\cdot Tr[\Upsilon{(\frac{\partial}{\partial
x}-\frac{\partial}{\partial y})}^2\Psi]=4\Psi''(|x-y|),$$ We
select $\Psi(|x-y|)$ such that
\begin{displaymath}
\Psi''(|x-y|)=\left\{\begin{array}{l}
\frac{\varphi}{3{|x-y|}^{1-\frac{2\sigma}{\ln E}}{[2\epsilon R_1
E^{-m\sigma(1+2\delta)}]}^{\frac{2\sigma}{\ln
E}}}{[(m+1)\sigma(1+2\delta)\ln E]}^{-1}, \\ \textrm{when
$0<|x-y|\le
2\epsilon R_1 E^{-m\sigma(1+2\delta)} $},\\
\frac{\varphi}{3|x-y|}{[\ln {\frac{2\epsilon R_1}{|x-y|}}
+\sigma(1+2\delta)\ln E]}^{-\frac 1 2}{[(m+1)\sigma(1+2\delta)\ln
E]}^{-\frac 1 2}, \\ \textrm{when $ 2\epsilon R_1
E^{-m\sigma(1+2\delta)} \le |x-y| \le 2\epsilon R_1 $},
\end{array} \right.
\end{displaymath}
Integrating the above expression we have
\begin{displaymath}
\Psi'(|x-y|)=\left\{\begin{array}{l}
\frac{\varphi}{3\frac{2\sigma}{\ln E}}{[\frac{|x-y|}{2\epsilon R_1
E^{-\sigma(1+2\delta)}}]}^{\frac{2\sigma}{\ln
E}}{[(m+1)\sigma(1+2\delta)\ln E]}^{-1},\\
\textrm{when $0<|x-y|\le
2\epsilon R_1 E^{-m\sigma(1+2\delta)} $},\\
\frac{\varphi}{3\frac{2\sigma}{\ln E}}{[(m+1)\sigma(1+2\delta)\ln
E]}^{-1}+2\\
-2{[\ln {\frac{2\epsilon R_1}{|x-y|}}+\sigma(1+2\delta)\ln
E]}^{\frac 1 2}{[(m+1)\sigma(1+2\delta)\ln E]}^{-\frac 1 2},\\
\textrm{when $ 2\epsilon R_1 E^{-m\sigma(1+2\delta)} \le |x-y| \le
2\epsilon R_1 $},
\end{array} \right.
\end{displaymath}
$$\Psi(|x-y|)=\int_0^{|x-y|}\Psi'(t)dt.$$
Hence (\ref{eq45}) is valid by the explicit expression of
$\Psi,\Psi'$.

We want to prove that
\begin{equation}\label{eq49}
-Tr(\Upsilon Z_2)\geq KJ^{\delta/2}\Psi''(|x-y|)
\end{equation}
under suitable conditions. (\ref{eq49}) follows from
$$J^{\delta/2-1}\frac{\sigma}{\ln E}\frac 1 {2\epsilon R_1}{(\frac{|x-y|}{2\epsilon
R_1})}^{\frac{2\sigma}{\ln E}-1}\ll J^{\delta/2}\Psi''(|x-y|),$$
or
$$I=\frac 3 \varphi \frac{\sigma}{\ln E}\frac 1 {2\epsilon R_1}{(\frac{|x-y|}{2\epsilon
R_1})}^{\frac{2\sigma}{\ln E}-1}\frac 1 {\Psi''(|x-y|)}\ll J,$$
$$={(\frac{|x+y-2x_1|}{2R_1})}^2+{(\frac{|x-y|}{2\epsilon
R_1})}^{\frac{2\sigma}{\ln E}}+e^{-2m(1+2\delta)\sigma}.$$ In the
interval
$$2\epsilon R_1 E^{-m\sigma(1+2\delta)} \le |x-y| \le
2\epsilon R_1,$$
$$I=\frac 3 \varphi \frac{\sigma}{\ln E}\frac 1 {2\epsilon R_1}{(\frac{|x-y|}{2\epsilon
R_1})}^{\frac{2\sigma}{\ln E}-1}|x-y|{[\ln {\frac{2\epsilon
R_1}{|x-y|}}+\sigma(1+2\delta)\ln E]}^{\frac 1 2}$$
$$\cdot{[(m+1)\sigma(1+2\delta)\ln E]}^{\frac 1 2}$$
$$\leq \frac 3 \varphi \sigma^2 (m+1)(1+2\delta){(\frac{|x-y|}{2\epsilon
R_1})}^{\frac{2\sigma}{\ln E}}\ll {(\frac{|x-y|}{2\epsilon
R_1})}^{\frac{2\sigma}{\ln E}},$$ if the following condition is
true
\begin{equation}\label{eq50}
(m+1)\sigma^2\ll 1.
\end{equation}
In the interval $0<|x-y|< 2\epsilon R_1 E^{-m\sigma(1+2\delta)}$,
$$I=\frac 3 \varphi \frac{\sigma}{\ln E}\frac 1 {2\epsilon R_1}{(\frac{|x-y|}{2\epsilon
R_1})}^{\frac{2\sigma}{\ln E}-1}{|x-y|}^{1-\frac{2\sigma}{\ln
E}}$$
$$\cdot{[2\epsilon R_1 E^{-m \sigma (1+2\delta)}]}^{\frac{2\sigma}{\ln
E}}[(m+1)\sigma(1+2\delta)\ln E]$$
$$=\frac 6 \varphi \sigma^2 (1+2\delta)E^{-m(1+2\delta)\frac{2\sigma}{\ln
E}}(m+1)$$
$$=\frac 6 \varphi (m+1) \sigma^2
(1+2\delta)e^{-2m(1+2\delta)\sigma}\ll e^{-2m(1+2\delta)\sigma},$$
if the condition (\ref{eq50}) is true. Hence, (\ref{eq49}) is
valid under the restriction (\ref{eq50}).

It is easy to prove that
$$\|Z_1\|\ll -Tr(\Upsilon Z_2)$$
when $$\sigma \epsilon^2 \ln E \ll 1 ;$$
$$\|Z_2\| \omega^2 \ll -Tr(\Upsilon Z_2)$$
when $$ \epsilon^{2\beta} \sigma \ln E {(\frac{|x-y|}{2\epsilon
R_1})}^{2\beta [1-\frac 1 2 \frac{\sigma}{\ln
E}]-\frac{2\sigma}{\ln E}}$$
$$\ll \epsilon^{2\beta} \sigma \ln E \ll 1 ;$$
$$\|Z\|\ll -Tr(\Upsilon Z_2)$$
when $$\epsilon^{\beta} \sigma \ln E \ll 1 ;$$
$$\omega \ll -Tr(\Upsilon Z_2)$$
when $$ \epsilon^{1+\beta} R_1 \sigma \ln E \ll 1 .$$ Taking
$\epsilon = {(\ln E)}^{-\frac 2 \beta}$, then the estimates for
$\|Z_1\|,\|Z_2\|,\|Z\|,\omega$ are all true. Hence we have the
following lemma.

\par {\bf Lemma 5.2} The estimates (\ref{eq46}) is valid when
(\ref{eq47}) and (\ref{eq50}) are true.
\begin{proof}
The lemma is proved since (\ref{eq23}) is true by the above
estimates.
\end{proof}
(\ref{eq46}) implies that
\begin{equation}\label{eq51}
\begin{aligned}
|u(x)-u(y)-<Du(x_1),x-y>|\leq \qquad \qquad \qquad \\
(1+\varphi)K{[{(\frac{|x+y-2x_1|}{2R_1})}^2
+{(\frac{|x-y|}{2\epsilon R_1})}^{\frac{2\sigma}{\ln
E}}+e^{-2m(1+2\delta)\sigma}]}^{\frac \delta 2}|x-y|,
\end{aligned}\end{equation} $\forall x,y \in Q.$ Taking constants
$m,\sigma$ such that $m\sigma$ is large and $m\sigma^2$ is small,
hence (\ref{eq50}) is valid. Taking $|x-y|$ small such that
\begin{equation}\label{eq52}
{(\frac{|x-y|}{2\epsilon R_1})}^{\frac{2\sigma}{\ln E}}\leq
e^{-2m(1+2\delta)\sigma}
\end{equation}
then (\ref{eq51}) implies that
\begin{equation}\label{eq53}
\begin{aligned}
|u(x)-u(y)-<Du(x_1),x-y>|\leq \qquad \qquad \\
(1+\varphi)K{[{(\frac{|x+y-2x_1|}{2R_1})}^2
+2e^{-2m(1+2\delta)\sigma}]}^{\frac \delta 2}|x-y|,
\end{aligned}\end{equation} $\forall x,y \in Q $ satisfies
(\ref{eq52}). We relax the restriction (\ref{eq52}) for $x,y$ and
assume $x,y$ satisfying
$$\forall x,y \in
\check{Q}:\{x,y|{(\frac{|x+y-2x_1|}{2R_1})}^2 +{(\frac{|x-y|}{2
R_1})}^{2} < 1-e^{-2m(1+2\sigma)\delta}\},$$ interpolating line
segment $\overline{xy}$ by $z_i=\frac{x+y} 2 + \frac{i}{2M}(x-y),$
where $i=-M, -M+1, \cdots, 0, \cdots, M-1, M$. By applying
(\ref{eq53}), we have
$$|u(z_{-i})-u(z_{-i-1})-<Du(x_1),z_{-i}-z_{-i-1}>|$$
$$+|u(z_{i+1})-u(z_i)-<Du(x_1),z_{i+1}-z_i>|$$
$$\leq
(1+\varphi)K[(U+V)^{\frac{\delta}{2}}+(U-V)^{\frac{\delta}{2}}]\frac{|x-y|}{2M},$$
where $$U={(\frac{|x+y-2x_1|}{2R_1})}^2
+{(\frac{2i+1}{2M})}^2{(\frac{|x-y|}{2R_1})}^2
+2e^{-2m(1+2\delta)\sigma},$$
$$V=\frac{2i+1}{M}<\frac{x+y-2x_1}{2R_1}, \frac{x-y}{2R_1}>.$$
Since $0<\frac{\delta}{2}\leq \frac{3}{4}\leq 1$, we have $$
{(U+V)}^{\frac{\delta}{2}}+{(U-V)}^{\frac{\delta}{2}}\leq
2U^{\frac{\delta}{2}}$$ $$\leq 2{[{(\frac{|x+y-2x_1|}{2R_1})}^2
+{(\frac{|x-y|}{2 R_1})}^{2} +
2e^{-2m(1+2\delta)}]}^{\frac{\delta}{2}}$$ Summing for
$i=0,1,\cdots,M-1$, we have
\begin{equation}\label{eq54}\begin{aligned}
|u(x)-u(y)-<Du(x_1),x-y>|\leq \qquad \qquad \\
(1+\varphi)K{[{(\frac{|x+y-2x_1|}{2R_1})}^2
+{(\frac{|x-y|}{2R_1})}^2+\xi^2]}^{\frac \delta 2}|x-y|,
\end{aligned}\end{equation}
$\forall x,y \in \tilde{Q}:\{x,y|{(\frac{|x+y-2x_1|}{2R_1})}^2
+{(\frac{|x-y|}{2R_1})}^2<1-\xi^2\}$ where
\begin{equation}\label{eq55}
2e^{-2m(1+2\delta)\sigma}\leq 2 e^{-3m\sigma}\equiv {\xi}^2 \ll 1.
\end{equation} $\forall \delta \in [\frac 1 2, \frac 3 2].$

Now we state and prove the fundamental lemma on decreasing of
Lipschitz coefficient in certain constant ratio accompanied with
decreasing of radius of spherical region in suitable constant
ratio.

\par {\bf Remark 5.1} Let $y\to x$ in (5.22), approximately we have

\begin{equation}\label{remark1}
|Du(x)-Du(x_1)|\leq
(1+\varphi+\varepsilon)(\frac{|x-x_1|}{R_1})^\delta),
(\varepsilon=\frac{\xi R_1}{|x-x_1|}),
\end{equation}
$\delta=\frac{3}{2}$ is a sharp estimate. But Eq.(\ref{remark1})
is valid only in the spherical shell $\xi R_1\leq|x-x_1|\leq R_1$,
and not for the sphere $|x-x_1|\leq R_1$, so that,
Eq.(\ref{remark1}) and (5.22) are reasonable.

\par {\bf Lemma 5.3} Let $x_1$ be a Caffarelli point of
$B(x_0,R_0)$ satisfying the Lipschitz condition (\ref{eq42}),
$\forall x,y \in B(x_1,\sqrt{2}R_1) \subset B(x_0,R_0)$.Under the
restriction (\ref{eq47}), $ \forall $ points $x_2 \in B(x_1,\xi
R_1)$, $\exists$ Caffarelli point $y_1$ of $B(x_1,R_1)$ such that
$y_1 \in B(x_1,\xi R_1)$, and $|y_1-x_2|\leq \frac{\xi}{3} R_1$,
we have
\begin{equation}\label{eq56}
|Du(y_1)-Du(x_1)|\leq \eta K,
\end{equation}
\begin{equation}\label{eq57}
|u(x)-u(y)-<Du(y_1),x-y>|\leq \zeta K|x-y|,
\end{equation}
$\forall x,y$ such that $\{x,y|{|\frac {x+y} 2 -x_2|}^2 +{(\frac
{|x-y|} 2)}^2\leq \frac {R_1^2} 9 \}$, where $\eta$ is a small
constant and $\zeta$ is a constant ,$\zeta <1$. Both $\eta,\zeta$
depend on $\xi$ only when we fixed $\delta$ to be $\delta=\frac 1
2$.

\begin{proof}
Since $x_2\in B(x_1,\xi R_1)$, in $B(x_2,\frac{\xi}{3} R_1)\cap
B(x_1,\xi R_1)$,  $\exists$ Caffarelli point of $B(x_1,R_1)$,
because of $B(x_2,\frac{\xi}{3} R_1)\cap B(x_1,\xi R_1) \supset
B(x_2+\frac{x_1-x_2}{|x_1-x_2|}\frac{\xi}{6}, \frac{\xi}{6}R_1)$
and the density of sphere
$B(x_2+\frac{x_1-x_2}{|x_1-x_2|}\frac{\xi}{6}, \frac{\xi}{6}R_1)$
with respect to $B(x_1,R_1)$ satisfying
  ${(\frac{\xi}{6})}^n \gg E^{-\Gamma}$ (If this condition
is not satisfied, since $\xi$ and $\Gamma$ are independent of $E$,
we substitute $E$ by a large one, such that this condition is
satisfied.)  We denote one of Caffarelli point by $y_1$. Hence
$y_1 \in B(x_1,\xi R_1), |y_1 -x_2|\leq \frac{\xi}{3} R_1$ and
moreover $Du(y_1)$ exists.
 Take $y=y_1$ in (\ref{eq54}) and then divide it by $|x-y_1|$
and let $x\rightarrow y_1$, we have
$$|Du(y_1)-Du(x_1)|\leq (1+\varphi){(2\xi^2)}^{\frac \delta 2} K $$
$$=(1+\varphi){(2\xi^2)}^{\frac 1 4} K\equiv \eta K ,$$ this is
(\ref{eq56}). Since
$$\{x,y|{|\frac {x+y} 2 -x_2|}^2 +{(\frac
{|x-y|} 2)}^2\leq \frac {R_1^2} 9 \} $$ $$\subset \{x,y|{|\frac
{x+y} 2 -x_1|}^2 +{(\frac {|x-y|} 2)}^2\leq \frac {R_1^2} 4 \}
\subset \tilde{Q}.$$ applying (\ref{eq54}) and (\ref{eq56}) we
have
$$|u(x)-u(y)-<Du(y_1),x-y>|\leq [\eta + (1+\varphi){(\frac 1 4
+\xi^2)}^{\frac 1 4}]K|x-y|\leq \zeta K|x-y|.$$
$$ \forall \{x,y|{|\frac {x+y} 2 -x_2|}^2 +{(\frac
{|x-y|} 2)}^2\leq \frac {R_1^2} 9 \} ,$$ where we denote
\begin{equation}\label{eq58}
\zeta = \eta +(1+\varphi)[{(\frac 1 4+   2\xi^2)}^{\frac 1
4}{(1+9\xi^2)}^{\frac 1 2}+\eta {(\frac 1 9 +2\xi^2)}^{-\frac 1
2}].\end{equation}
 We take a little larger $\zeta$ over our
necessary for the sake of next lemma. Since $\xi, \eta$ are small,
we have
$$\zeta \thickapprox (1+\varphi){(\frac 1 4)}^{\frac 1 4}=[1+\frac
1 2 (\sqrt{2}-1)]\frac{ 1 }{\sqrt{2}}=\frac 1 2 (1+\frac {1}
{\sqrt{2}})<1,$$ hence we can take $\xi$ small such that
$\zeta<1$.

The lemma is proved completely.
\end{proof}

\par {\bf Remark 5.2} The only reason for approximating the
general point $x_2$ by Caffarelli point $x_1,y_1,\ldots$ is that
there exists first derivative on Caffarelli points.

\par {\bf Lemma 5.4} We have $u \in C^{1+\alpha}(B(x_0,\frac{\xi}{\sqrt{2}} R_0))$,
where
\begin{equation}\label{eq59}
\alpha=\frac {\ln \zeta}{\ln \frac 1 3}\end{equation}
\begin{proof}
Fixed a $x_2 \in (B(x_0,\frac{\xi}{\sqrt{2}}R_0))$. Take a
Caffarelli point $x_1$ of $(B(x_0,\frac{\xi}{\sqrt{2}} R_0))$ such
that $|x_1-x_2|\leq \frac{\xi}{\sqrt{2}} R_0$. Take
$R_1=(1-\xi)\frac{1}{\sqrt{2}} R_0$ and $K=2M_1$. Denote
$R^{(0)}=R_1, R^{(k)}=\frac {R_1}{3^k}, K^{(0)}=K, K^{(k)}=\zeta^k
K (k=0,1,2,\cdots)$. Denote $x_1\equiv y_0$, take Caffarelli point
$y_1$ of $B(y_0,R^{(0)})$ such that $y_1 \in B(y_0,\xi
R^{(0)})\cap B(x_2,\frac{\xi}{3} R^{(0)})$, hence $|y_1-x_2|<
\frac{\xi}{3} R^{(0)}=\xi R^{(1)},$ i.e. $y_1\in B(x_2,\xi
R^{(1)})$, take Caffarelli point $y_2$ of $B(y_1,R^{(1)})$ such
that $y_2 \in B(y_1,\xi R^{(1)})\cap B(x_2,\frac{\xi}{3} R^{(1)})
\cdots $ In general, take Caffarelli point $y_k$ of
$B(y_{k-1},R^{(k-1)})$ such that $y_k \in B(y_{k-1},\xi
R^{(k-1)})\cap B(x_2,\frac{\xi}{3} R^{(k-1)})$, $\forall
k=1,2,\cdots.$

Since we restrict $R_0$ to be $\leq M_1$, we have
$$R^{(k)}=\frac {R_1}{3^k}\leq \frac {R_0}{3^k \sqrt{2}}\leq \frac {M_1}{3^k}\leq
K \zeta^k = K^{(k)} (k=0,1,2,\cdots),$$  i.e. (\ref{eq47}) is
valid for all $k=0,1,2,\cdots$.

We prove by induction that
\begin{equation}\label{eq60}
|Du(y_k)-Du(y_{k-1})|\leq K^{(k)},
\end{equation}
\begin{equation}\label{eq61}
|u(x)-u(y)-<Du(y_k),x-y>|\leq K^{(k)}|x-y|,
\end{equation}
$\forall x,y$ such that $\{ {|\frac {x+y} 2 -x_2|}^2 +{(\frac
{|x-y|} 2)}^2\leq {(R^{(k)})}^2 \}$, $\forall k=1,2,\cdots$

When $k=1$, (\ref{eq60}) follows from (\ref{eq56}) and $\eta K\leq
\zeta K=K^{(1)}$, (\ref{eq61}) follows from (\ref{eq57}).

When (\ref{eq60}) and (\ref{eq61}) are valid for $k-1$,
substituting $y_0=x_1, R^{(0)}=R_1, K^{(0)}=K=2M_1$ by $y_{k-1},
R^{(k-1)}, K^{(k-1)}$ in lemma 5.2, 5.3. Since $R^{(k)}\leq
K^{(k)}$, hence lemma 5.2 is valid and lemma 5.3 is valid by using
$y_k, R^{(k)}$ and $K^{(k)}$ instead of $y_1, R^{(1)}=\frac {R_1}
3, K^{(1)}=\zeta K$. (\ref{eq60}) and (\ref{eq61}) follow for $k$
by applying lemma 5.3 and $\eta K^{(k-1)}\leq \zeta
K^{(k-1)}=K^{(k)}$. i.e. (\ref{eq60}) and (\ref{eq61}) are valid
$\forall k=1,2,\cdots$.

Applying (\ref{eq60}) we have
$$ \sum_{k=1}^\infty |Du(y_k)-Du(y_{k-1})|\leq \sum_{k=1}^\infty K
\zeta ^k <\infty.$$ Hence the series
$$\sum_{k=1}^\infty [Du(y_k)-Du(y_{k-1})] $$ converges, denote its
limit by $\tilde{a}-Du(x_1)$. $\forall x,y$ satisfies $\{x,y
{|\frac {x+y} 2 -x_2|}^2 +{(\frac {|x-y|} 2)}^2\leq {(R^{(k)})}^2
\}$, applying (\ref{eq61}) we have
\begin{equation}\label{eq62} \begin{aligned}
|u(x)-u(y)-<\tilde{a},x-y>| \leq |u(x)-u(y)-<Du(y^k),x-y>| \\
+ |<Du(y^k)-\tilde{a},x-y>| \qquad\qquad\qquad\\
\leq K^{(k)}|x-y|+(\sum_{l=k+1}^\infty K\zeta^l)|x-y|=\frac
{1+\zeta}{1-\zeta}K^{(k)}|x-y| \\
=\frac {1+\zeta}{1-\zeta}K \zeta^k |x-y|=\frac {1+\zeta}{1-\zeta}K
{(\frac{R^{(k)}}{R^{(0)}})}^\alpha |x-y|\qquad\\
\leq \frac {1+\zeta}{1-\zeta}K
{(\frac{3R^{(k+1)}}{(1-\xi)R^{(0)}})}^\alpha |x-y| \leq CR^\alpha
|x-y|,\qquad
\end{aligned}\end{equation}
when we denote
$$R= {({|\frac {x+y} 2 -x_2|}^2 +{(\frac
{|x-y|} 2)}^2)}^{\frac 1 2},$$
\begin{equation}\label{eq63}
C=\frac {1+\zeta}{1-\zeta} 2 M_1 {(\frac {3}{1-\xi})}^\alpha
{R_0}^{-\alpha},\end{equation}
and assume
$$ R^{(k+1)}\leq R\leq R^{(k)} (k=0,1,2,\cdots)$$

Put the case $k=0,1,2,\cdots$ together, we obtain that
(\ref{eq62}) is true
$$\forall 0\leq {({|\frac {x+y} 2 -x_2|}^2 +{(\frac
{|x-y|} 2)}^2)}^{\frac 1 2}\leq R\leq R^{(1)}=\frac {R_1} 3=
\frac{1-\xi}{3\sqrt{2}} R_0.$$

Putting $y=x_2$ and dividing (\ref{eq62}) by $|x-x_2|$, then let
$x\rightarrow x_2$, we have $Du(x_2)$ exists and
$\tilde{a}=Du(x_2)$. Hence we have
\begin{equation}\label{eq64}
|u(x)-u(y)-<Du(x_2),x-y>|\leq CR^\alpha |x-y|, \end{equation}
$\forall x,y$ satisfies $\{x,y|{({|\frac {x+y} 2 -x_2|}^2 +{(\frac
{|x-y|} 2)}^2)}^{\frac 1 2}\leq R\leq R_1 =\frac{1-\xi}{3\sqrt{2}}
R_0\}$

(\ref{eq64}) implies that
$$|Du(x)-Du(x_2)|\leq C{|x-x_2|}^\alpha, \forall x,x_2 \in
B(x_0,\frac{\xi}{\sqrt{2}} R_0)$$ Hence
$$u \in C^{1+\alpha}(B(x_0,\frac{\xi}{\sqrt{2}} R_0))$$
The lemma is proved.
\end{proof}

Now we study the case that Lipschitz coefficient contains a factor
$R^\gamma$, where $\gamma \in (0,1)$.

\par {\bf Lemma 5.5}  $\forall x_2 \in B(x_0,\frac{\xi}{\sqrt{2}} R_0)$, if constant
$\gamma \in (0,1)$ and positive constant $H,S\leq
\frac{(1-\xi)}{3\sqrt{2}}R_0$ exist such that we have
\begin{equation}\label{eq65}
|u(x)-u(y)-<Du(x_2),x-y>|\leq H R^\gamma |x-y|,\end{equation}
$\forall x,y$ satisfies $\{x,y|{({|\frac {x+y} 2 -x_2|}^2 +{(\frac
{|x-y|} 2)}^2)}^{\frac 1 2}\leq R\leq S\}$

Then in case \begin{equation}\label{eq66} \gamma \leq 1-\alpha,
\end{equation} where $\alpha$ is defined by (\ref{eq59}).
$\exists$ constant $\tilde{H}(>H)$ and $\tilde{S}(<S)$ such that
we have
\begin{equation}\label{eq67}
|u(x)-u(y)-<Du(x_2),x-y>|\leq H R^{\gamma+\alpha} |x-y|,
\end{equation}
$\forall x,y$ satisfies $\{x,y|{({|\frac {x+y} 2 -x_2|}^2 +{(\frac
{|x-y|} 2)}^2)}^{\frac 1 2}\leq R\leq \tilde{S}\}.$

\begin{proof}
Denote

$R^{(0)}=min\{S,H^{\frac 1 {1-\gamma}}\},
K^{(0)}=H{R^{(0)}}^\gamma,$ where the meaning of $R^{(0)}$ is
different from that in lemma 5.4.

Let $y_0 \in B(x_2,\xi R^{(0)})$ be a Caffarelli point of $u(x)$
in $B(x_2,R^{(0)})$. In the region $B(y_0,(1-\xi )R^{(0)})$, the
lemma 5.2 for estimate $|u(x)-u(y)-<Du(y_0),x-y>|$ is also valid
in the present case, this is because of (\ref{eq47})
$$R^{(0)} \leq H{(R^{(0)})}^\gamma =K^{(0)} $$ is true. Take
$K=H{(R^{(0)})}^\gamma, \delta=\frac 1 2 +\gamma $. From
(\ref{eq54}),

$\forall x,y \in \tilde{Q} :\{x,y|{(\frac {|x+y-2y_0|}
{2R^{(0)}})}^2 +{(\frac {|x-y|} {2R^{(0)}})}^2 \leq 1-\xi^2 \},$

we have
$$|u(x)-u(y)-<Du(y_0),x-y>|$$ $$\leq (1+\varphi) H{(R^{(0)})}^\gamma
{[{(\frac {|x+y-2y_0|} {2R^{(0)}})}^2 +{(\frac {|x-y|}
{2R^{(0)}})}^2+2\xi^2]}^{\frac 1 4 +\frac \gamma 2}|x-y|.$$

Hence we have for all Caffarelli point $y_1$ of $B(y_0,R^{(0)})$
satisfying $y_1 \in B(y_0,\xi R^{(0)})\cap B(x_2,\frac{\xi}{3}
R^{(0)})$, we have
\begin{equation}\label{eq68}
|Du(y_1)-Du(y_0)|\leq (1+\varphi) H{(R^{(0)})}^\gamma
{(2\xi^2)}^{\frac 1 4 +\frac \gamma 2} \leq (1+\varphi) K^{(0)}
{(2\xi^2)}^{\frac 1 4 +\frac \gamma 2}=\eta K^{(0)},
\end{equation}
$\forall x,y$ satisfying $\{x,y|{|\frac {x+y} 2 -x_2|}^2 +{(\frac
{|x-y|} 2)}^2 <{(\frac {R^{(0)}} 3)}^2 \}$, we have
\begin{equation}\label{eq69}\begin{aligned}
|u(x)-u(y)-<Du(y_1),x-y>| \qquad\qquad \\ \leq
[(1+\varphi)H{(R^{(0)})}^\gamma {(\frac 1 9 +2\xi^2)}^{\frac
\gamma
2} {(\frac 1 4 +2\xi^2)}^{\frac 1 4}+\eta K^{(0)}]|x-y|\\
\leq \zeta \frac {K^{(0)}}{3^\gamma} |x-y|, \qquad\qquad\qquad
\end{aligned}\end{equation} where $\zeta$ is defined by (\ref{eq58}).

(\ref{eq68}) and (\ref{eq69}) are the similar relation of lemma
5.3 in the present case.

Define $ R^{(k)}=\frac{R^{(0)}}{3^k}, K^{(k)}=H{(R^{(k)})}^\gamma
\zeta^k (k=1,2,3,\cdots)$

Take Caffarelli point $y_k$ of $B(y_{k-1},R^{(k-1)})$ such that
$y_k \in B(y_{k-1},\xi R^{(k-1)})\cap B(x_2,\frac{\xi}{3}
R^{(k-1)})$ .

First we verify (\ref{eq47}) in the present case. Applying
(\ref{eq66}) we have
$$ 3^{1-\gamma}\zeta \geq 3^\alpha \zeta=3^{\frac{\ln{\frac 1
\zeta}}{\ln 3}}\zeta =\frac{1}{\zeta} \zeta =1 ,$$ or it is
$$\frac 1 {3^k} \leq \frac{\zeta^k}{{(3^k)}^\gamma}.$$  Hence
$$R^{(k)}=\frac{R^{(0)}}{3^k}\leq
\frac{K^{(0)}}{3^k}=\frac{H{(R^{(0)})}^\gamma}{3^k}$$
$$\leq H{(\frac{R^{(0)}}{3^k})}^\gamma \zeta^k=H{(R^{(k)})}^\gamma\zeta^k=K^{(k)},$$
i.e. (\ref{eq47}) is valid in the present case.

Then we prove (\ref{eq60}) and (\ref{eq61}) by induction. When
$k=1$, applying (\ref{eq68}) we have
$$|Du(y_1)-Du(y_0)|\leq \eta K^{(0)} =\frac{\eta}{\zeta}3^\gamma
K^{(1)} \leq \frac{3\eta}{\zeta}K^{(1)}\leq K^{(1)}.$$ Applying
(\ref{eq69}) we have
$$|u(x)-u(y)-<Du(y_1),x-y>|\leq \zeta \frac {K^{(0)}}{3^\gamma} |x-y|
=K^{(1)}|x-y|,$$ $\forall x,y$ satisfies $\{x,y|{|\frac {x+y} 2
-x_2|}^2 +{(\frac {|x-y|} 2)}^2 <{(\frac {R^{(0)}} 3)}^2
={(R^{(1)})}^2 \}$

 i.e. (\ref{eq60}) and (\ref{eq61}) are valid
for $k=1$.

Suppose (\ref{eq60}) and (\ref{eq61}) are valid for $k-1$. Since
$K^{(k)}=H{(R^{(k)})}^\gamma \zeta^k$, take $\delta=\frac 1 2
+\gamma$ and apply (\ref{eq54}),
\begin{equation}\label{eq70}\begin{aligned}
|u(x)-u(y)-<Du(y_k),x-y>| \qquad\qquad \\
\leq (1+\varphi)K^{(k)}
 {[{(\frac {|x+y-2y_k|} {2R^{(k)}})}^2 +{(\frac
{|x-y|} {2R^{(k)}})}^2+2\xi^2]}^{\frac 1 4 +\frac \gamma 2}|x-y|,
\end{aligned}\end{equation}
$\forall x,y$ satisfies $$\{x,y|{(\frac {|x+y-2y_k|}{2R^{(k)}})}^2
+{(\frac {|x-y|} {2R^{(k)}})}^2 < 1-\xi^2 \}.$$ Hence applying
(\ref{eq70}) we have
$$|Du(y_{k+1})-Du(y_k)|\leq (1+\varphi)K^{(k)}{(2\xi^2)}^{\frac 1
4 +\gamma}\leq (1+\varphi)K^{(k)}{(2\xi^2)}^{\frac 1 4}$$
$$=\eta K^{(k)} =\frac{\eta}{\zeta}3^\gamma
K^{(k+1)} \leq \frac{3\eta}{\zeta}K^{(k+1)}\leq K^{(k+1)}.$$ And
$$|u(x)-u(y)-<Du(y_{k+1}),x-y>| $$
$$ \leq (1+\varphi) [K^{(k)}{(\frac 1 9
+2\xi^2)}^{\frac \gamma 2} {(\frac 1 4 +2\xi^2)}^{\frac 1 4}+\eta
K^{(k)}]|x-y|$$
$$=\frac{\zeta}{3^\gamma}K^{(k)} |x-y|=K^{(k+1)}|x-y|,$$
$\forall x,y$ satisfies $\{x,y|{|\frac {x+y} 2 -x_2|}^2 +{(\frac
{|x-y|} 2)}^2 <{(\frac {R^{(k)}} 3)}^2 ={(R^{(k+1)})}^2 \}.$

Hence (\ref{eq60}) and (\ref{eq61}) are valid for all
$k=1,2,\cdots$.

Since $Du(x) \in C, \forall x \in B(x_0,\frac{\xi}{\sqrt{2}}
R_0)$. Applying (\ref{eq60}) and (\ref{eq61}) we have
$$|u(x)-u(y)-<Du(x_2),x-y>| $$
$$\leq |u(x)-u(y)-<Du(y_k),x-y>|+[\sum_{l=k+1}^\infty
K^{(l)}]|x-y|$$
$$\leq K^{(k)}|x-y| +[\sum_{l=k+1}^\infty
K^{(l)}]|x-y|$$
$$=H[{(R^{(k)})}^\gamma
\zeta^k+\sum_{l=k+1}^\infty {(R^{(l)})}^\gamma \zeta^l]|x-y|$$
$$=\frac{1+\frac{\zeta}{3^\gamma}}{1-\frac{\zeta}{3^\gamma}}H {(R^{(k)})}^\gamma
\zeta^k |x-y|$$
$$=\frac{1+\frac{\zeta}{3^\gamma}}{1-\frac{\zeta}{3^\gamma}}H {(R^{(k)})}^\gamma
{(\frac{R^{(k)}}{R^{(0)}})}^\alpha |x-y|.$$ In the region
$R^{(k+1)}\leq {({|x-x_2|}^2+{|y-x_2|}^2)}^{\frac 1 2}=R\leq
R^{(k)}$, we have
$$R^{(k)}=3 R^{(k+1)}\leq 3 R.$$
Hence
$$|u(x)-u(y)-<Du(x_2),x-y>|\leq \tilde{H} R^{\gamma+\alpha}|x-y|,$$
where
\begin{equation}\label{eq71}
\tilde{H}=\frac{1+\frac{\zeta}{3^\gamma}}{1-\frac{\zeta}{3^\gamma}}3^{\gamma+\alpha}
H{(R^{(0)})}^{-\alpha}. \end{equation}
 Putting $k=1,2,3,\cdots$
together we have
$$|u(x)-u(y)-<Du(x_2),x-y>|\leq \tilde{H} R^{\gamma+\alpha}|x-y|,$$
$\forall x,y$ satisfies $\{x,y|{|x-x_2|}^2+{|y-x_2|}^2\leq R\leq
\check{S}\},$
\begin{equation}\label{eq72}
\check{S}=\frac{R^{(0)}}{3} \leq \frac{1-\xi}{3}
\min\{S,H^{\frac{1}{1-\gamma}}\}. \end{equation}

The lemma is proved.
\end{proof}

\par {\bf Lemma 5.6} $\forall (x_0,R_0)\subset \Omega$, we have
$$u(x) \in C^{1+1}(B(x_0,\frac{\xi}{\sqrt{2}} R_0)),$$ where $\xi$
is defined by (\ref{eq55}).

\begin{proof}
Since $\alpha$ and $\zeta$ are defined by (\ref{eq59}) and
(\ref{eq58}), substituting $\zeta$ by a little large one such that
$\zeta<1$ still valid and $\frac 1 \alpha $ is a positive integer
$M$. Applying lemma 5.4 once we obtain (\ref{eq64}). Then applying
lemma 5.5 successively.

Denote
\begin{equation}\label{eq73}
R_j^{(0)}=\min \{S_j,H_j^{\frac 1 {1-j\alpha}}\}, \qquad
j=1,2,\cdots,M.
\end{equation}
where $S_j,H_j$ are the value of $S,H$ in lemma 5.5 corresponding
to $\gamma=j \alpha$.

Applying lemma 5.4 we have
$$S_1=\frac{1-\xi}{3\sqrt{2}} R_0,$$
$$H_1=\frac{1+\zeta}{1-\zeta} 2 M_1 {(\frac 3 {1-\xi})}^\alpha
{(\frac{1-\xi}{\sqrt{2}}R_0)}^{-\alpha}.$$ Since $M_1\geq R_0$, it
is easy to obtain that
\begin{equation}\label{eq74}
H_1\geq S_1^{1-\alpha}.
\end{equation}
Hence applying (\ref{eq73}) we have
$$ R_1^{(0)}= S_1.$$
Applying lemma 5.5 we have
\begin{equation}\label{eq75}
S_{j+1}=\frac{1-\xi}{3} \min\{S_j,H_j^{\frac 1 {1-j \alpha}}\}.
\end{equation}
\begin{equation}\label{eq76}
H_{j+1}=\frac{1+\frac{\zeta}{3^{j \alpha}}}{1-\frac{\zeta}{3^{j
\alpha}}} 3^{(j+1)\alpha}H_j{(R_j^{(0)})}^{-\alpha} \geq
H_j{(R_j^{(0)})}^{-\alpha}.
\end{equation}
We prove by induction that
\begin{equation}\label{eq77}
H_j\geq S_j^{1-j \alpha}
\end{equation}
is true, then by (\ref{eq73}) we have $R_j^{(0)}=S_j$.

When $j=1$, (\ref{eq77}) is true by (\ref{eq74}). If (\ref{eq77})
is true for $j$, applying (\ref{eq76}) we have
$$ H_{j+1}\geq H_j{(R_j^{(0)})}^{-\alpha}\geq S_j^{1-j \alpha}
S_j^{-\alpha} = S_j^{1-(j+1)\alpha},$$ i.e. (\ref{eq77}) is valid
for $j$ substitute by $j+1$.

Hence (\ref{eq77}) is true $\forall j=1,2,\cdots, M$. Applying
(\ref{eq75}) we have $$ S_{j+1} =\frac{1-\xi}{3} S_j
(j=1,2,\cdots,M-1).$$ Hence $$ S_{j+1}={(\frac{1-\xi}{3})}^j
S_1={(\frac{1-\xi}{3})}^{j+1}\frac{R_0}{\sqrt{2}}.$$
Especially we
have
$$S_M=R_M^{(0)}={(\frac{1-\xi} 3)}^M \frac{R_0}{\sqrt{2}}=[1+o(1)]e^{-\frac{\ln
3}{\alpha}}\frac{R_0}{\sqrt{2}} $$
$$\thickapprox e^{ \frac{{(\ln 3)}^2}{\ln [{\frac 1 2}(1+\frac 1
{\sqrt{2}})]}}\frac{R_0}{\sqrt{2}},$$ when $\xi$ is small. $e^{
\frac{{(\ln 3)}^2}{\ln [{\frac 1 2}(1+\frac 1 {\sqrt{2}})]}}$ is a
constant $<1$ and is independent of $\xi$. Hence
$$|u(x)-u(y)-<Du(x_2),x-y>|\leq \tilde{C} \frac{R}{R_0}|x-y|$$
$$\forall x,y \in B(x_0,R_M^{(0)}).$$ Since $R_M^{(0)} > \xi
\frac{R_0}{\sqrt{2}},$ hence we have $\forall x_2,x_3 \in
B(x_0,\frac{\xi}{\sqrt{2}} R_0),$
$$|Du(x_3)-Du(x_2)|\leq \tilde{C} \frac{|x_3-x_2|}{R_0}$$

The lemma is proved.
\end{proof}

Our main theorem 2.1.

\begin{proof}
Theorem 2.1 follows by applying theorem 5.1 and lemma 5.7.
\end{proof}

\par {\bf A final remark} After we proved the regularity result $u\in
C^{1+1}(\Omega)$, a natural question has arisen. Is the regularity
result the best possible or not? If it is yes, we must give an
example to show that the solution of (\ref{eq1}) under the
conditions (\ref{eq2})-(\ref{eq4})) cannot belong to
$C^2(\Omega)$. More exactly, the solution cannot possesses more
regularity when it is lack of concavity assumption, so that more
regularity assumptions on $F$ has no effect. Hence we present the
following problem.
\par {\bf An open problem}
Let the equation (\ref{eq1}) satisfying the following assumptions:

(a) $F(x,z,p,X)$ is sufficiently smooth with respect to its
arguments $x,z,p,X$.

(b) Assumption (\ref{eq2}).

(c) $F$ is strict monotone with respect to $z$ in the following
meaning
$$F(x,z,p,X)+K(w-z)\leq F(x,w,p,X),$$

(d) $F$ satisfies the natural structure condition
$$|F|+(1+|p|)|F_p|+|F_z|+\frac 1{1+|p|}|F_x|$$
$$\leq \mu (1+|p|^2+|X|),$$
where $\lambda , \Lambda, K, \mu$ are positive constants.

 Can you construct an example to
show that the solution $u(x)$ of (\ref{eq1}) under the assumptions
(a)-(d) does not belong to $C^2(\Omega)$?

\bigskip

\institution[ {Department of Mathematics, Zhejiang University,
Hangzhou, China(310027)}

%[2]{\\Tongji University, Shanghai, China} [3]{\\Department of
%Mathematics, Zhejiang University, Hangzhou, China(310027)}
\email{dogc@zju.edu.cn}

%% references info

\onecolumn

\pagebreak
\end{document}